\documentclass[amssymb,11pt,reqno]{amsart}

\newcommand{\R}{\mathbb{R}}
\newcommand{\C}{\mathbb{C}}

\newcommand{\N}{\mathbb{N}}
\newcommand{\G}{\Gamma}
\newcommand{\rth}{R_{\theta}}
\newcommand{\test}{\mathcal S}
\newcommand {\D}{\Delta}

\newtheorem*{main}{Main Theorem}
\newtheorem{thm}{Theorem}

\newtheorem{lemma}{Lemma}
\newtheorem{prop}{Proposition}

\usepackage{graphicx}
\usepackage{latexsym}

\def\vol{\mbox{\rm Vol}}
\def\r2n{\R^{2n}}

\begin{document}

\vspace{1.5 cm}

%\begin{center}{ Last modified Dec 2, 2004  } \end{center}

%%%%%%%%%%%%%%%%%%%%%%%%%%%%%%%%%%%%%%%%%%%%%%%%%%%%%%%%%%%%%%%%%%%%
\title[The Complex Busemann-Petty problem for arbitrary measures]
      {The Complex Busemann-Petty problem for arbitrary measures.}

\author{Marisa Zymonopoulou}
\address{
Marisa Zymonopoulou\\
Department of Mathematics\\
University of Missouri\\
Columbia, MO 65211, USA}

\email{marisa@@math.missouri.edu}
%%%%%%%%%%%%%%%%%%%%%%%%%%%%%%%%%%%%%%%%%%%%%%%%%%%%%%%%%%%%%%%%%%%%\begin{abstract}
\begin{abstract}The complex Busemann-Petty problem asks whether origin symmetric convex bodies in $\C^n$ with smaller central hyperplane sections necessarily have smaller volume. The answer is affirmative if $n\leq 3$ and negative if $n \geq 4.$ In this article we show that the answer remains the same if the volume is replaced by an ``almost'' arbitrary measure. This result is the complex analogue of Zvavitch's generalization to arbitrary measures of the original real Busemann-Petty problem.

\end{abstract}

\maketitle

\section{Introduction}

In 1956 the Busemann-Petty problem was posed (see [BP]), asking the following question:
suppose that $K$ and $L$ are two origin symmetric convex bodies in $\R^n$ such that for every $\xi \in S^{n-1},$

$$\vol_{n-1}\bigl(K\cap \xi^{\perp}\bigr) \leq \vol_{n-1}\bigl(L\cap \xi^{\perp}\bigr).$$
Does it follow that
$$\vol_{n}\bigl(K\bigr) \leq \vol_{n}\bigl(L \bigr) \ \ ?$$

\noindent
The answer is affirmative if $n\leq 4$ and negative if $n\geq 5.$ The problem was solved in the late 90's as a result of a series of papers ([LR], [Ba], [Gi], [Bu], [Lu], [Pa], [Ga], [Zh1], [K1], [K2], [Zh2], [GKS]; see [K5, p.3] for the history of the solution).

A few years later Zvavitch [Zv] showed that one can replace the volume by essentially any measure on $\R^n.$ Namely, if we
consider any even continuous positive function $f$ on $\R^n$ and denote by $\mu$ the measure with density $f,$ we can define
\begin{center}$\mu(D)=\int_{D}f(x)dx$ \hspace{10pt} and \hspace{10pt} $\mu(D\cap \xi^{\perp})=\int_{D\cap \xi^{\perp}}f(x)dx,$
\end{center}
for every closed bounded invariant with respect to all $\rth$ set $D$ in $\R^n$ and every $\xi \in S^{n-1}.$ 
 Then the Busemann-Petty problem for general measures is stated as follows:

 Suppose that $K$ and $L$ are two origin symmetric convex bodies in $\R^n$ such that, for every $\xi \in S^{n-1},$
$$\mu(K\cap \xi^{\perp}) \leq \mu(L\cap \xi^{\perp}).$$

\noindent
Does it follow that $$\mu(K)\leq \mu(L) \ ?$$
%the integral of $f$ over any central hyperplane section of $K$ is smaller than the same
%for $L.$ Does it imply that the integral of $f$ over $K$ is smaller than that over $L?$
Surprisingly, the answer remains the same as in the original problem. It is affirmative for $n\leq 4$ and negative for $n\geq 5.$

Zvavitch's ideas for general measures were applied and further developed in [R], [Y1] and [Y2], for hyperbolic and spherical spaces and for sections of lower dimensions.

In this article we study the complex version of the Busemann-Petty problem for arbitrary measures. 

Let $\xi \in \C^n$ with $ |\xi|=1.$ We denote by
$$H_{\xi}=\{z\in \C^n \ : \ (z,\xi)=\sum\limits_{k=1}^{n}z_k \overline{\xi}_k=0\}$$

\noindent
the complex hyperplane perpendicular to $\xi.$

Origin symmetric convex bodies in $\C^n$ are the unit balls of norms on $\C^n.$ We denote by $\|\cdot\|_K$ the norm corresponding to the body $K$
$$K=\{z \in \C^n \ : \|z\|_{K}\leq 1\}.$$

\noindent
We identify $\C^n$ with $\R^{2n}$ using the mapping
$$\xi=(\xi_1,\ldots ,\xi_n)=(\xi_{11}+i\xi_{12},\ldots ,\xi_{n1}+i\xi_{n2})\longmapsto (\xi_{11},\xi_{12},\ldots, \xi_{n1},\xi_{n2})$$
\noindent
and observe that under this mapping the complex hyperplane $H_{\xi}$ turns into a $(2n-2)$-dimensional subspace of $\R^{2n}$ orthogonal to the vectors

   \begin{center}
   $\xi=(\xi_{11},\xi_{12},\ldots, \xi_{n1},\xi_{n2})$  \hspace{1.1pt} and \hspace{1.1pt}  $\xi^{\perp}=(-\xi_{12},\xi_{11},\ldots, -\xi_{n2},\xi_{n1}).$
   \end{center}

\noindent
Since norms on $\C^n$ satisfy the equality
$$\|\lambda z\|=|\lambda|\|z\|, \ \ \forall z \in \C^n, \ \forall \lambda \in \C^n,$$
origin symmetric complex convex bodies correspond to those origin symmetric convex bodies $K$ in $\R^{2n}$ that are invariant with respect to any coordinate-wise two-dimensional rotation, namely for each $\theta \in [0,2\pi]$ and each $x=(x_{11},x_{12},\ldots, x_{n1},x_{n2}) \in \R^{2n}$
\begin{equation}
\|x\|_{K}=\|\rth(x_{11},x_{12}), \ldots, \rth(x_{n1},x_{n2})\|_{K},
\end{equation}

\noindent
where $\rth$ stands for the counterclockwise rotation of $\R^2$ by the angle $\theta$ with respect to the origin. If a convex body satisfies $(1)$ we will say that \emph{it is invariant with respect to all $\rth$}.

The complex Busemann-Petty problem ([KKZ]) can now be formulated as follows:
Suppose $K$ and $L$ are origin symmetric invariant with respect to all $\rth$
 convex bodies in $\R^{2n}$ such that
 $$\vol_{2n-2}(K\cap H_\xi)\leq \vol_{2n-2}(L\cap H_\xi)$$
 for each $\xi$ from the unit sphere $S^{2n-1}$ of $\R^{2n}.$ Does it follow that
 $$\vol_{2n}(K) \leq \vol_{2n}(L) \  ?$$
 
As it is proved in [KKZ] the answer is affirmative if $n\leq 3$ and negative if $n\geq 4.$ 

Let $f$ be an even positive and continuous function on $\r2n$. We define a measure $\mu$ on $\r2n$ with density $f,$ so that
\begin{center}$\mu(D)=\int_{D}f(x)dx$ \hspace{10pt} and \hspace{10pt} $\mu(D\cap H)=\int_{D\cap H}f(x)dx$
\end{center}
for every closed bounded invariant with respect to all $\rth$ set $D$ in $\r2n$ and every $(2n-2)$-dimensional subspace $H$ of $\r2n.$ 
%As in the complex Busemann-Petty problem, the property of invariance of the complex convex bodies, when viewed as bodies in $\R^{2n},$ plays an important role to the solution of the problem. 
As it is proved in Section 3 (Lemma \ref{finv}), one may assume, without loss of generality, that the density $f$ is also invariant with respect to all rotations $\rth.$ We will call such a function \emph{$\rth$-invariant}.
Then, the complex Busemann-Petty problem for arbitrary measures is stated as follows:

\noindent
Suppose $K$ and $L$ are origin symmetric invariant with respect to all $\rth$ convex bodies in $\R^{2n}$ so that for every $\xi \in S^{2n-1}$
$$\mu(K\cap H_{\xi})\leq \mu(L\cap H_{\xi}),$$
does it follow that
$$\mu(K)\leq\mu(L) \ ?$$

In this article we prove that, analogously to the real case, the solution remains the same for arbitrary measures with a positive continuous density. 

 Note that, the positivity assumption on $f$ is necessary, because otherwise one may assume that the density is identically zero where the affirmative answer to the problem holds trivially in all dimensions.
% Namely, if $f$ is the even continuous density function of a measure $\mu$ on $\R^{2n},$ we can define its average over the unit circle

%\begin{equation}\label{eqt:ave}
%\tilde{f}(x)=\frac{1}{2\pi}\int_0^{2\pi}f(\rth x)d\theta, \ \forall x \in \R^{2n}.\end{equation}

%\noindent
%Then $\tilde{f}$ is an even continuous and invariant with respect to all $\rth$ function on $\R^{2n}$ such that

%$$\mu(D)=\int_{D}\tilde{f}(x)dx $$

%\noindent
%for every closed bounded invariant with respect to all $\rth$ set $D$ in $\R^{2n}.$

\medskip

\section{The Fourier analytic connection to the problem}

Through out this paper we use the Fourier transform of distributions. The Schwartz class of rapidly decreasing infinitely differentiable functions (test functions) in $\R^n$ is denoted by $\test(\R^n),$ and the space of distributions over $\test(\R^n)$ by $\test^{\prime}(\R^n).$ The Fourier transform $\hat{f}$ of a distribution $f \in \test^{\prime}(\R^n)$ is defined by $\langle\hat{f},\phi\rangle=\langle f,\hat{\phi}\rangle$ for every test function $\phi.$ A distribution is called even homogeneous of degree $p \in \R$ if $\langle f(x),\phi(x/\alpha)\rangle=|\alpha|^{n+p}\langle f,\phi\rangle$ for every test function $\phi$ and every $\alpha \in \R, \ \alpha\neq 0.$ The Fourier transform of an even homogeneous distribution of degree $p$ is an even homogeneous distribution of degree $-n-p.$ A distribution $f$ is called positive definite if, for every test function $\phi, \ \langle f, \phi\ast\overline{\phi(-x)}\rangle \geq 0.$ By Schwartz's generalization of Bochner's theorem, this is equivalent to $\hat{f}$ being a positive distribution in the sense that  $\langle\hat{f},\phi\rangle \geq 0$ for every non-negative test function $\phi,$ (see [K5, section 2.5] for more details).

 A compact set $K \subset\R^n$ is called a star body, if every straight line that passes through the origin crosses the boundary of the set at exactly two points and the boundary of $K$ is continuous in the sense that the \emph{Minkowski functional} of $K,$ defined by
$$\|x\|_K=\min \{\alpha \geq 0 : x \in \alpha K \}$$
is a continuous function on $\R^n.$

A star body $K$ in $\R^n$ is called $k$-smooth (infinitely sooth) if the restriction of $\|x\|_{K}$ to the sphere $S^{n-1}$ belongs to the class of $C^k(S^{n-1} )\  (C^{\infty}(S^{n-1})).$ It is well-known that one can approximate any convex body in $\R^n$ in the radial metric,
$d(K, L)=\sup \{|\rho_{K}(\xi)-\rho_{L}(\xi)|,\ \xi \in S^{n-1} \},$
by a sequence of infinitely smooth convex bodies. The proof is based on a simple convolution argument (see for example [Sch, Theorem 3.3.1]). It is also easy to see that any convex body in $\R^{2n}$ invariant with respect to all $\rth$ rotations can be approximated in the radial metric by a sequence of infinitely smooth convex bodies invariant with respect to all $\rth.$ This follows from the same convolution argument, because invariance with respect to $\rth$ is preserved under convolutions.

If $D$ is an infinitely smooth origin symmetric star body in $\R^n$ and $0<k <n,$ then the Fourier transform of the distribution $\|x\|_D^{-k}$ is a homogeneous function of degree $-n+k$ on $\R^n,$ whose restriction to the sphere is infinitely smooth (see [K5, Lemma 3.16]).

The following Proposition is a spherical version of Parseval's formula established in [K3], (see also [K5, Lemma 3.22]):
\begin{prop}\label{prop:parseval}
Let $D$ be an infinitely smooth origin symmetric star body in $\R^n$ and $g \in C^{k-1}(\R^n)$ even homogeneous of degree $-n+k$ function. Then
$$\int_{S^{n-1}}g(\theta)\| \theta\|_{D}^{-k}d\theta=(2\pi)^n
\int_{S^{n-1}}\hat{g}(\xi)\bigl(\| \theta\|_{D}^{-k}\bigr)^{\wedge}(\xi) d\xi.$$
\end{prop}

\medskip

The concept of an intersection body was introduced by Lutwak [Lu]. This concept was generalized in [K3], as follows:
Let $1\leq k <n,$ and let $D$ and $L$ be two origin symmetric star bodies in $\R^n.$ We say that $D$ is the \emph{$k$-intersection body of $L$} if for every $(n-k)-$dimensional subspace $H$ of $\R^n$
$$\vol_k(D\cap H^{\perp})=\vol_{n-k}(L\cap H ).$$
We introduce the class of $k$-intersection bodies, as those star bodies that can be obtained as the limit, in the radial metric, of a sequence of $k$-intersection bodies of star bodies. 
 A Fourier analytic characterization of $k$-intersection bodies was proved in [K4].

%More generally, we say that an origin symmetric star body $D\subset \R^n$ is a $k$-intersection body if there exists a finite Borel measure $\mu$ on $S^{n-1}$ so that for every test function $\phi \in \test (\R^n),$
%$$\int_{\R^n}\|x\|_D^{-k}\phi(x)dx=
%\int_{S^{n-1}}\Bigl(\int_0^{\infty}t^{k-1}\hat{\phi}(t\xi)dt\Bigr)d\mu(\xi).$$

\noindent

\begin{prop}\label{prop:k-posdef}
An origin symmetric star body $D$ in $\R^n$ is a $k$-intersection body, $1\leq k \leq n-1,$ if and only if $\|\cdot\|_{D}^{-k}$ is a positive definite distribution.

\end{prop}

\medskip

\bigskip

Let $1\leq k<2n$ and let $H$ be an $(2n-k)-$dimensional subspace of $\R^{2n}.$ We denote by $\chi(\cdot)$ the indicator function on $[-1,1]$ and by $|\cdot|_2$ the Euclidean norm in the proper space. We fix an orthonormal basis $e_1,\ldots,e_k$ in the orthogonal subspace $H^{\perp}.$ For any convex body $D$ in $\R^{2n}$ and any even positive continuous function $f$ on $\R^{2n}$ we define the $(2n-k)-$dimensional parallel section function $A_{f,D,H}$ as a function on $\R^k$ such that

\begin{equation}\label{eqt:secf}A_{f,D,H}(u)=
\int_{\{x\in \R^{2n}:(x,e_1)=u_1,\ldots,(x,e_k)=u_k\} }\chi(\|x\|_D) f(x)dx, \ u\in \R^{k}. \end{equation}

The original lower dimensional parallel section function that corresponds to the $(n-k)$-dimensional volume  of the section of $D$ with a subspace $H$ (put $n$ instead of $2n$ and $f=1$), was defined in [K4].
Note that at $0$ the function $A_{f,D,H}$ measures the central section of the body $D$ by the subspace $H.$ Passing to polar coordinates on $H$ we have that 

\begin{eqnarray}\label{eqn:polarH}
A_{f,D,H}(0)&=&\mu(D\cap H)=\int_{H}\chi(\|x\|_D) f(x)dx \nonumber \\
            &=&\int\limits_{S^{2n-1}\cap H}\Bigl(\int_{0}^{\|\theta\|_D^{-1}}r^{2n-3}f(r\theta)dr\Bigr)d\theta.
            \end{eqnarray}

If $D$ is infinitely smooth and $f \in \C^{\infty}(\R^{2n}),$ the function $A_{f,D,H}$ is infinitely differentiable at the origin (see [K5, Lemma 2.4]). So we can consider the action of the distribution $|u|_2^{-q-k}/ \G(-q/2)$ on $A_{f,D,H}$ and apply a standard regularization argument (see for example [K5, p.36] and [GS, p.10]). Then the function

\begin{equation}\label{eqt:q-A}
q \longmapsto \left< \frac{|u|_2^{-q-k}}{\G(-\frac{q}{2})},A_{f,D,H}(u)\right>
\end{equation}

\noindent
is an entire function of $q\in \C.$
If $q=2m, \ m \in \N\cup \{0\},$ then
$$\left<\frac{|u|_2^{-q-k}}{\G(-\frac{q}{2})}\Big|_{q=2m}, A_{f,D,H}(u)\right>$$
\begin{equation*}\label{eqt:q=2m}=\frac{(-1)^m |S^{k-1}|}{2^{m+1}k(k+2)\cdots(k+2m-2)}\Delta^m A_{f,D,H}(0),
\end{equation*}
where $|S^{k-1}|=2\pi^{k/2}/\G(k/2)$ is the surface area of the unit sphere in $\R^k,$ and $\Delta=\sum_{i=1}^k \partial^2/\partial u_i^2$ is the $k$-dimensional Laplace operator (see [GS, p.71-74]).
Note that the function (\ref{eqt:q-A}) is equal, up to a constant, to the fractional power of $\Delta^{q/2}A_{f,D,H}$ (see [KKZ] or [K4] for complete definition).

\medskip

\noindent
\textbf{Remark.} If a body $D$ is $m$-smooth (or infinitely smooth) and $f \in C^m(\R^{2n})$ (or $C^{\infty}(\R^{2n}))$ it is easy for one to see that the function
$$x \mapsto |x|_2^{-m}\int_0^{\frac{|x|_2}{\|x\|_K}}r^{2n-3}f\bigl(r\frac{x}{|x|}\bigr)dr$$ is also $m$-times (infinitely) continuously differentiable on $\R^{2n}\setminus \{0 \}.$

\medskip

The proof of following proposition is similar to that of Proposition 4 in [KKZ].
So we omit it here.

%We use a well-known formula (see [GS, p.76]): for every $v \in \R^k$ and $q<-k+1,$
%$$(v_1^2+\cdots +v_k^2)^{\frac{-q-k}{2}}$$
%\begin{equation}\label{formulaGS}
%=\frac{\G(-q/2)}{2\G((-q-k+1)/2)\pi^{(k-1)/2}}\int_{S^{k-1}}|(v,u)|^{-q-k}du.
%\end{equation}

\begin{prop}\label{prop:A}
Let $D$ be an infinitely smooth origin symmetric convex body in $\R^{2n}, \ f \in C^{\infty}(\R^{2n}),$ and $1\leq k < 2n.$  Then for every $(2n-k)-$dimensional subspace $H$ of $\R^{2n}$ and any $q \in \R, \ -k<q< 2n-k,$

$$\left< \frac{|u|_2^{-q-k}}{\G(-\frac{q}{2})}, A_{f,D,H}(u)\right>$$
\begin{equation}\label{eqt:A^{q}}=\frac{2^{-q-k}\pi^{-\frac{k}{2}}}{\G\bigl(\frac{q+k}{2}\bigr)}
\int\limits_{S^{2n-1}\cap H^{\perp}}
\Bigl(|x|_2^{-2n+k+q}\int_{0}^{\frac{|x|_2}{\|x\|_D}}
r^{2n-k-1-q}f\bigl(r\frac{x}{|x|_2}\bigr)dr\Bigr)^{\wedge}(\theta)d\theta.
\end{equation}

\noindent
Now, if $m \in \N \cup \{0\},$

$$\D^m A_{f,D,H}(0)$$
\begin{equation}\label{eqt:D^m}=\frac{(-1)^m}{(2\pi)^k}\int\limits_{S^{2n-1}\cap H^{\perp}}
\Bigl(|x|_2^{-2n+k+2m}\int_{0}^{\frac{|x|_2}{\|x\|_D}}
r^{2n-k-1-2m}f\bigl(r\frac{x}{|x|_2}\bigr)dr\Bigr)^{\wedge}(\theta)d\theta
\end{equation}
\end{prop}

\bigskip
The following (elementary) inequality is similar to Lemma 1 in [Zv]. 

\begin{lemma}\label{lm:elem}
Let $a, b>0$ and let $\alpha$ be a non-negative function on $(0, \max \{a,b \} ]$ so that the integrals below converge. Then
\begin{equation}\label{eqt:elem}
\int_{0}^{a}t^{2n-1}\alpha(t)dt-a^2\int^{a}_0t^{2n-3}\alpha(t)dt\leq \int_{0}^{b}t^{2n-1}\alpha(t)dt-a^2\int^{b}_0t^{2n-3}\alpha(t)dt.\!\! \end{equation}
\end{lemma}

%\bigskip

\section{Connection with $k$-intersection bodies}\label{measuresections}

 As mentioned in the Introduction, we can assume that the density function is $\rth$-invariant. This simple observation plays an important role to the solution of the problem. 

\begin{lemma}\label{finv}
Suppose $f$ is an even non-negative continuous function on $\R^{2n}$ and $\mu$ is a measure with density $f.$ Then there exists an even non-negative continuous function $\tilde{f}$ that is invariant with respect to all rotations $\rth$ such that
\begin{center}
$\mu(D)=\int_{D}\tilde{f}(x)dx $ \hspace{1.1pt}    and \hspace{1.1pt}    $\mu(D\cap H_{\xi})=\int_{D\cap H_{\xi}}\tilde{f}(x)dx,$
\end{center}
for every closed bounded invariant with respect to all $\rth$ set $D$ in $\R^{2n}$ and $\xi \in S^{2n-1}.$

\end{lemma}

\noindent
\textbf{Proof.}
We define its average over the unit circle, $\tilde{f}(x)=\frac{1}{2\pi}\int_0^{2\pi}f(\rth x)d\theta,$ for every $x\in \R^{2n}.$ Then
for every compact invariant with respect to all $\rth$ set $D$ in $\R^{2n},$
\begin{eqnarray*}\int_{D}\tilde{f}(x)dx&=&\frac{1}{2\pi}\int_{D}\int_0^{2\pi}f(\rth x)d\theta dx \\
&=&\frac{1}{2\pi}\int_0^{2\pi}\int_{\rth^{-1}D}f(y)dyd\theta=
\mu(D)\end{eqnarray*}
since $\rth^{-1}D=D,$ for all $\theta \in [0,2\pi].$

Moreover, since central sections of complex convex bodies by complex hyperplanes correspond to convex bodies in $\R^{2n-2}$ that are also invariant with respect to the $\rth$ rotations, we similarly get that for every $\xi \in S^{2n-1},$
$$\mu(D\cap H_{\xi})=\int_{D\cap H_{\xi}}\tilde{f}(x)dx.$$
\qed

Now, we are ready to express the measure of the central sections in terms of the Fourier transform.

\begin{thm}\label{thm:musection}
Suppose $K$ is an infinitely smooth origin symmetric invariant with respect to all $\rth$ convex body in $\R^{2n}, \ n \geq 2,$ and $f$ is an infinitely differentiable even positive and $\rth$-invariant function on $\R^{2n}.$ Then for every $\xi \in S^{2n-1}$

\begin{equation}\label{eqt:musection}
\mu(K\cap H_{\xi})=\frac{1}{2\pi}\Bigl(|x|_2^{-2n+2}\int_{0}^{\frac{|x|_2}{\|x\|_{K}}}
r^{2n-3}f\bigl(r\frac{x}{|x|_2}\bigr)dr\Bigr)^{\wedge}(\xi)
\end{equation}
\end{thm}

\noindent
In order to prove Theorem \ref{thm:musection} we need the following:

\begin{lemma}\label{lemma:constf}
Let $K$ and $f$ as in Theorem \ref{thm:musection}. Then for every $\xi \in S^{2n-1}$ the Fourier transform of the distribution
\begin{equation}\label{invdistr}|x|_2^{-2n+2}\int_0^{\frac{|x|_2}{\|x\|_{K}}} r^{2n-3}f (r\frac{x}{|x|_{2}})dr
\end{equation} is a constant function on $S^{2n-1}\cap H_{\xi}^{\perp}.$
\end{lemma}

\noindent
\textbf{Proof.}
The function $\|x\|^{-1}_{K}$ is invariant with respect to all $\rth$ (see Introduction), so, since $f$ is $\rth$-invariant it is easy to see that the distribution in (\ref{invdistr}) is a continuous function which is also invariant with respect to all rotations $\rth.$ By the connection between the Fourier transform of distributions and linear transformations, its Fourier transform is also invariant with respect to all $\rth$. As mentioned in the Introduction, the space $H_{\xi}^{\perp}$ is spanned by the vectors $\xi$ and $\xi^{\perp}.$ So every vector in $S^{2n-1}\cap H_{\xi}^{\perp}$ is a rotation $\rth,$ for some $\theta \in [0,2\pi],$ of $\xi$ and hence the Fourier transform of $$|x|_2^{-2n+2}\int_0^{\frac{|x|_2}{\|x\|_{K}}} r^{2n-3}f (r\frac{x}{|x|_{2}})dr$$ is a constant function on $S^{2n-1}\cap H_{\xi}^{\perp}.$
\qed

\medskip

\noindent
\textbf{Proof of Theorem \ref{thm:musection}}.
Let $\xi \in S^{2n-1}.$ In formula (\ref{eqt:D^m}) we put $ \  H_{\xi}=H, \ k=2$ and $ m=0.$ Then, by the definition of the lower dimensional section function $A_{f,D,H}(0),$ equation (\ref{eqn:polarH}), we have that

$$\mu(K\cap H_{\xi})=\frac{1}{(2\pi)^2}\int_{S^{2n-1}\cap H_{\xi}^{\perp}}\Bigl(|x|_2^{-2n+2}\int_{0}^{\frac{|x|_2}{\|x\|_{K}}}
r^{2n-3}f\bigl(r\frac{x}{|x|_2}\bigr)dr\Bigr)^{\wedge}(\eta)d\eta.$$

By Lemma \ref{lemma:constf}, the function under the integral is constant on the circle $S^{2n-1}\cap H_{\xi}^{\perp}.$ Since $\xi \in H_{\xi}^{\perp}$ we have that
$$\mu(K\cap H_{\xi})=\frac{1}{(2\pi)^2}2\pi\Bigl(|x|_2^{-2n+2}\int_{0}^{\frac{|x|_2}{\|x\|_{K}}}
r^{2n-3}f\bigl(r\frac{x}{|x|_2}\bigr)dr\Bigr)^{\wedge}(\xi)$$

\noindent
which proves the theorem.\qed

%The following (elementary) lemma, similar to lemma 1 in [Zv],(see also [K5, Lemma 5.14], will play an important role to the positive answer of the problem. Here we state the lemma without the proof.

%\begin{lemma}\label{lm:elem}
%Let $a, b>0$ and let $\alpha$ be a non-negative function on $(0, \max{a,b}]$ so that the integrals below converge. Then
%\begin{equation}\label{eqt:elem}
%\int_{0}^{a}t^{2n-1}\alpha(t)dt-a^2\int^{a}_0t^{2n-3}\alpha(t)dt\leq %\int_{0}^{b}t^{2n-1}\alpha(t)dt-a^2\int^{b}_0t^{2n-3}\alpha(t)dt\end{equation}
%\end{lemma}

\bigskip

As in the case of the complex Busemann-Petty problem the property of a body to be a $2$-intersection body is closely related to the solution of the complex Busemann-Petty problem for arbitrary measures.

\begin{thm}\label{thm:main2}
The solution of the complex Busemann-Petty problem for arbitrary measures in $\C^n $ has an affirmative answer if and only if every origin symmetric invariant with respect to all $\rth$ convex body in $\R^{2n}$ is a $2$-intersection body.
\end{thm}

The proof of Theorem \ref{thm:main2} will follow from the Remarks and the next lemmas.

\medskip

\noindent
\textbf{Remark 1.} To prove the affirmative part of the problem it is enough to consider infinitely smooth origin symmetric invariant with respect to all $\rth$ bodies. This is true because one can approximate, in the radial metric, from inside the body $K$  and from outside the body $L$ by infinitely smooth convex invariant with respect to all $\rth$ bodies. Then if the affirmative answer holds for infinitely smooth bodies it also holds in the general case.

\medskip

\noindent
\textbf{Remark 2.} Let $D$ be an origin symmetric convex body which is not a $k$-intersection body. Then, there exists a sequence of infinitely smooth convex bodies with strictly positive curvature which are not $k$-intersection bodies that converges in the radial metric to $D,$ (see [K5, Lemma 4.10]).  If, in addition, $D$ is invariant with respect to all $\rth,$ one can choose a sequence of bodies with the same property.

\medskip

\noindent
\textbf{Remark 3.} A simple approximation argument allows us to prove Theorem \ref{thm:main2} only for measures whose density is an infinitely differentiable even positive and $\rth$-invariant function on $\R^{2n}$. Let $f$ be the even positive continuous $\rth$-invariant density function of a measure $\mu,$ as it is defined in the Introduction. Then there exists an increasing sequence $g_n$ of even positive functions in $C^{\infty}(\R^{2n})$ such that  $g_n(x)\chi(\|x\|_D) \rightarrow f(x)\chi(\|x\|_D),$ a.e., for every compact set $D.$ Then by the Monotone Convergence Theorem we have that

\begin{center}$\int_{\R^{2n}} g_n(x)\chi(\|x\|_D)dx\rightarrow \mu(D)$ \hspace{1.1pt}
and \hspace{1.1pt} $\int_{H} g_n(x)\chi(\|x\|_D)dx\rightarrow \mu(H\cap D),$
\end{center} 

\noindent
as $n\rightarrow \infty,$ for every subspace $H$ of $\R^{2n}.$ In addition, by Lemma \ref{finv}, we may assume that every $g_n$ is also $\rth$-invariant.

\medskip
Now we are ready to prove the affirmative part of the complex Busemann-Petty problem for arbitrary measures.

\begin{lemma}\label{lm:pos}
Suppose $K$ and $L$ are infinitely smooth origin symmetric invariant with respect to all $\rth$ convex bodies in $\R^{2n}$ so that $K$ is a $2$-intersection body and let $f$ be an infinitely differentiable even positive $\rth$-invariant function on $\R^{2n}.$ Then, if for every $\xi \in S^{2n-1}$
\begin{equation}\label{eqt:sectineq}
\mu(K\cap H_{\xi}) \leq \mu(L\cap H_{\xi})
\end{equation}
then
$$\mu(K)\leq \mu(L).$$
\end{lemma}

\noindent
\textbf{Proof.}
By the remark before Proposition \ref{prop:A} and [K5, Lemma 3.16], the Fourier transform of the distributions
\begin{center}
$|x|_2^{-2n+2}\int_{0}^{\frac{|x|_2}{\|x\|_{K}}}
r^{2n-3}f\bigl(r\frac{x}{|x|_2}\bigr)dr,$ \hspace {7pt} and \hspace{7pt} $|x|_2^{-2n+2}\int_{0}^{\frac{|x|_2}{\|x\|_{L}}}
r^{2n-3}f\bigl(r\frac{x}{|x|_2}\bigr)dr$
\end{center}
%and $\|x\|_K^{-2}$
are homogeneous of degree $-2$ and continuous functions on $\R^{2n}\setminus \{0\}.$ So, by Theorem \ref{thm:musection}, the inequality (\ref{eqt:sectineq}) becomes

$$\Bigl(|x|_2^{-2n+2}\int_{0}^{\frac{|x|_2}{\|x\|_{K}}}
r^{2n-3}f\bigl(r\frac{x}{|x|_2}\bigr)dr\Bigr)^{\wedge}(\xi)$$
$$\leq \Bigl(|x|_2^{-2n+2}\int_{0}^{\frac{|x|_2}{\|x\|_{L}}}
r^{2n-3}f\bigl(r\frac{x}{|x|_2}\bigr)dr\Bigr)^{\wedge}(\xi).$$
Since $K$ is an infinitely smooth $2$-intersection body, by Proposition \ref{prop:k-posdef} and [K5, Theorem 3.16] the Fourier transform of the distribution $\|x\|_K^{-2}$ is a non-negative continuous, outside the origin, function on $\R^{2n}.$ Multiplying both sides of the latter inequality by $\bigl(\|x\|_K^{-2}\bigr)^{\wedge}$ and applying the spherical version of Parseval, Proposition \ref{prop:parseval}, we have that

$$\int_{S^{2n-1}}\bigl(\|x\|^{-2}_{K}\bigr)^{\wedge}(\xi)\Bigl(|x|_2^{-2n+2}\int_{0}^{\frac{|x|_2}{\|x\|_{K}}}
r^{2n-3}f\bigl(r\frac{x}{|x|_2}\bigr)dr\Bigr)^{\wedge}(\xi)d\xi$$
$$\leq \int_{S^{2n-1}}\bigl(\|x\|^{-2}_{K}\bigr)^{\wedge}(\xi)\Bigl(|x|_2^{-2n+2}\int_{0}^{\frac{|x|_2}{\|x\|_{L}}}
r^{2n-3}f\bigl(r\frac{x}{|x|_2}\bigr)dr\Bigr)^{\wedge}(\xi)d\xi,$$
which gives

$$\int_{S^{2n-1}}\|x\|^{-2}_{K}\int_{0}^{\|x\|_{K}^{-1}}
r^{2n-3}f(rx)drdx$$
\begin{equation}\label{eqt:secintineq}\leq \int_{S^{2n-1}}\|x\|^{-2}_{K}\int_{0}^{\|x\|_{L}^{-1}}
r^{2n-3}f(rx)drdx.\end{equation}

\noindent
We use the elementary inequality, equation (\ref{eqt:elem}), with $a=\|x\|_{K}^{-1}, b=\|x\|_{L}^{-1}$ and $\alpha(r)=f(rx)$ and integrate over $S^{2n-1}.$ Then
$$\int_{S^{2n-1}}\bigl(\int_{0}^{\|x\|_{K}^{-1}}r^{2n-1}f(rx)dr\bigr)dx-
\int_{S^{2n-1}}\|x\|_{K}^{-2}\bigl(\int_{0}^{\|x\|_{K}^{-1}}r^{2n-3}f(rx)dr\bigr)dx$$
\begin{equation}\label{eqt:sphelem}
\leq \int_{S^{2n-1}}\bigl(\int_{0}^{\|x\|_{L}^{-1}}r^{2n-1}f(rx)dr\bigr)dx-
\int_{S^{2n-1}}\|x\|_{K}^{-2}\bigl(\int_{0}^{\|x\|_{L}^{-1}}r^{2n-3}f(rx)dr\bigr)dx
\end{equation}

\noindent
We add the equations (\ref{eqt:secintineq}) and (\ref{eqt:sphelem}) and have that
$$\int_{S^{2n-1}}\bigl(\int_0^{\|x\|_{K}^{-1}}r^{2n-1}f(rx)dr\bigr)dx \leq \int_{S^{2n-1}}\bigl(\int_0^{\|x\|_{K}^{-1}}r^{2n-1}f(rx)dr\bigr)dx$$
which immediately implies that
$$\mu(K)\leq \mu(L).$$
\qed

For the negative part we need a perturbation argument to construct a body that will give a counter-example to the problem. The following lemma (without the assumption of invariance with respect to $\rth$ rotations) was proved in [Zv, Proposition 2] (see also [K5, Lemma 5.16]). The new body immediately inherits the additional property of invariance with respect to all $\rth$ of the original convex body.

\begin{lemma}\label{lemma:countexmpl}
Let $L$ be an infinitely smooth origin symmetric convex body with positive curvature  and let $f, g \in C^2(\R^{2n}),$  such that $f$ is strictly positive on $\R^{2n}.$ For $\varepsilon >0$ we define a star body $K$ so that
$$\int_0^{\|x\|^{-1}_K} t^{2n-3}f(tx)dt=\int_0^{\|x\|_{L}^{-1}}t^{2n-3}f(tx)dt-\varepsilon g(x), \
\forall x \in S^{2n-1}.$$
Then, if $\varepsilon$ is small enough the body $K$ is convex. Moreover, if $L$ is invariant with respect to all $\rth,$ and $f, g$ are $\rth$-invariant then $K$ is also invariant with respect to all $\rth.$

\end{lemma}

\medskip

\begin{lemma}\label{lm:neg}
Let $f \in C^{\infty}(\R^{2n})$ is an even positive $\rth$-invariant function. Suppose $L$ is an infinitely smooth origin symmetric invariant with respect to all $\rth$ convex body in $\R^{2n}$ with positive curvature which is not a $2$-intersection body. Then there exists an origin symmetric invariant with respect to all $\rth$ convex body $K$ in $\R^{2n}$ so that for every $\xi \in S^{2n-1}$
$$\mu(K\cap H_{\xi}) \leq \mu(L\cap H_{\xi})$$
but
$$\mu(K)> \mu(L).$$
\end{lemma}

\textbf{Proof.}
The body $L$ is infinitely smooth, so, by [K5, Lemma 3.16],  the Fourier transform of $\|x\|_{L}^{-2}$ is a continuous function on $\R^{2n}.$ Since $L$ is not a $2$-intersection body, by Proposition \ref{prop:k-posdef} there exists an open set $\Omega \subset S^{2n-1}$ where the Fourier transform of $\|x\|_L^{-2}$ is negative. We can assume that $\Omega$ is invariant with respect to rotations $\rth$ since $L$ is.

 Using a standard perturbation procedure for convex bodies, see for example [KKZ, Lemma 5] and [K5, p.96], we define an even non-negative invariant with respect to all $\rth$ function
$h \in C^{\infty}(S^{2n-1})$ whose support is in $\Omega.$ We extend $h$ to an even homogeneous function $h(\frac{x}{|x|_2})|x|_2^{-2}$ of degree $-2$ on $\R^{2n}.$ Then, by [K5, Lemma 3.16] the Fourier transform of $h(\frac{x}{|x|_2})|x|_2^{-2}$ is an even homogeneous function $g(\frac{x}{|x|_2})|x|_2^{-2n+2}$ of degree $-2n+2$ on $\R^{2n},$ with $g\in C^{\infty}(S^{2n-1}).$ Moreover, $g$ is also invariant with respect to rotations $\rth.$ 

 The assumptions for the body $L$ allow us to apply Lemma \ref{lemma:countexmpl} and take $\varepsilon >0$ small enough to define a convex body $K$ by
$$|x|_2^{-2n+2}\int_0^{\frac{|x|_2}{\|x\|_K}}t^{2n-3}f\bigl(t\frac{x}{|x|_2}\bigr)dt$$
$$=|x|_2^{-2n+2}\int_0^{\frac{|x|_2}{\|x\|_L}}t^{2n-3}f\bigl(t\frac{x}{|x|_2}\bigr)dt-
\varepsilon g\bigl(\frac{x}{|x|_2}\bigr)|x|_2^{-2n+2}.$$
We apply Fourier transform to both sides of the latter inequality. Then, by Theorem \ref{thm:musection}, since $h\geq 0,$ we obtain the following inequality for the measures of the central sections of $K$ and $L$ by the subspace $H_{\xi},$
$$\mu(K\cap H_{\xi})=\frac{1}{2\pi}\Bigl(|x|_2^{-2n+2}\int_{0}^{\frac{|x|_2}{\|x\|_{K}}}
r^{2n-3}f\bigl(r\frac{x}{|x|_2}\bigr)dr\Bigr)^{\wedge}(\xi)$$
$$=\frac{1}{2\pi}\Bigl(|x|_2^{-2n+2}\int_{0}^{\frac{|x|_2}{\|x\|_{L}}}
r^{2n-3}f\bigl(r\frac{x}{|x|_2}\bigr)dr\Bigr)^{\wedge}(\xi)-(2\pi)^{2n-1}\varepsilon h(\xi)$$
$$\leq \mu(L\cap H_{\xi})$$
On the other hand, the function $h$ is positive only where $\bigl(\|\cdot\|^{-2}_{L}\bigr)^{\wedge}$ is negative. So, for every $\xi \in S^{2n-1},$
$$\bigl(\|\cdot\|^{-2}_{L}\bigr)^{\wedge}(\xi)
\Bigl(|x|_2^{-2n+2}\int_{0}^{\frac{|x|_2}{\|x\|_{K}}}
r^{2n-3}f\bigl(r\frac{x}{|x|_2}\bigr)dr\Bigr)^{\wedge}(\xi)$$
$$=\bigl(\|\cdot\|^{-2}_{L}\bigr)^{\wedge}(\xi)\Bigl(|x|_2^{-2n+2}\int_{0}^{\frac{|x|_2}{\|x\|_{L}}}
r^{2n-3}f\bigl(r\frac{x}{|x|_2}\bigr)dr\Bigr)^{\wedge}(\xi)$$
$$-(2\pi)^{2n}\bigl(\|\cdot\|^{-2}_{L}\bigr)^{\wedge}(\xi)\varepsilon h(\xi)$$
$$>
\bigl(\|\cdot\|^{-2}_{L}\bigr)^{\wedge}(\xi)\Bigl(|x|_2^{-2n+2}\int_{0}^{\frac{|x|_2}{\|x\|_{L}}}
r^{2n-3}f\bigl(r\frac{x}{|x|_2}\bigr)dr\Bigr)^{\wedge}(\xi),$$
 Now, we integrate the latter inequality over $S^{2n-1}$ and apply the spherical version of Parseval's identity. Then similarly to Lemma \ref{lm:pos}, we apply the elementary inequality for integrals, Lemma \ref{lm:elem}, and conclude that
$$ \mu(K)> \mu(L).$$
\qed

\bigskip

\section{The solution of the problem }\label{solCBPGM}
To prove the main result of this paper we need to determine the dimensions in which an origin symmetric invariant with respect to all $\rth$ convex body in $\R^{2n}$ is a $2$-intersection body.

\begin{main}
The solution to the complex Busemann-Petty problem for arbitrary measures is affirmative if $n\leq 3$ and negative if $n\geq 4.$
\end{main}

\textbf{Proof.} It is known that an origin symmetric invariant with respect to $\rth,$ convex body in $\R^{2n}, n\geq 2,$ is a $k$-intersection body if $k\geq 2n-4$ (see [KKZ]). Hence, we obtain an affirmative answer to the complex Busemann-Petty problem for arbitrary measures if $n\leq 3.$

Now, suppose that $n \geq 4.$ The unit ball $B_q^n$ of the complex space $l_q^n, \ q>2,$ considered as a subset of $\R^{2n}:$
$$B^n_q=\{x\in \R^{2n}:\|x\|_q=\bigl( (x_{11}^2+x_{12}^2)^{q/2}+\cdots +(x_{n1}^2+x_{n2}^2)^{q/2}\bigr)^{1/q} \leq 1 \}$$
provides a counter-example for the Lebegue measure ($f=1$), of a body that is not a $k$-intersection body for $k<2n-4$ (see [KKZ, Theorem 4]). By Proposition \ref{prop:k-posdef} this implies that for $n\geq 4$ the distribution $\|x\|_q^{-2}$ is not positive definite. Then the result follows by Theorem \ref{thm:main2}.
\qed

\bigbreak

{\bf Acknowledgments:}
The author was partially supported by the
NSF grant DMS-0652571. Also, part of this work was carried out when the author was visiting the Pacific Institute of Mathematics, which the author thanks
for its hospitality.

\end{document}